\documentclass[aos,MSNbibl,nameyear,dvips]{arximspdf}
\usepackage{graphicx}

\doi{10.1214/10-AOS853}
\volume{39}
\issue{2}
\pubyear{2011}
\firstpage{702}
\lastpage{730}

\makeatletter

\newtheorem{Theorem}{Theorem}[section]
\newtheorem{Lemma}[Theorem]{Lemma}

\newproclaim{Remark}[Theorem]{Remark}
\newproclaim{Example}[Theorem]{Example}

\newcommand{\bdoi}{\doi}

\newcommand{\Leb}{\operatorname{Leb}}

\newcommand{\bs}{\mathbf}

\newcommand{\FF}{\mathcal{F}}
\newcommand{\HH}{\mathcal{H}}
\newcommand{\MM}{\mathcal{M}}
\newcommand{\QQ}{\mathcal{Q}}
\newcommand{\SSS}{\mathcal{S}}
\newcommand{\XX}{\mathcal{X}}

\newcommand{\R}{\mathbb{R}}

\newcommand{\Ex}{\mathbb{E}}

\newcommand{\argmax}{\mathop{\arg\max}}

\newcommand{\csupp}{\operatorname{csupp}}
\newcommand{\conv}{\operatorname{conv}}
\newcommand{\interior}{\operatorname{interior}}
\newcommand{\dom}{\operatorname{dom}}

\newcommand{\eps}{\varepsilon}

\makeatother

\begin{document}
\begin{frontmatter}

\title{Approximation by log-concave distributions, with applications to regression}
\runtitle{Log-concave approximations}

\begin{aug}
\author[A]{\fnms{Lutz} \snm{D\"{u}mbgen}\corref{}\thanksref{t1}\ead[label=e1]{duembgen@stat.unibe.ch}},
\author[B]{\fnms{Richard} \snm{Samworth}\ead[label=e2]{r.samworth@statslab.cam.ac.uk}}\\
and
\author[A]{\fnms{Dominic} \snm{Schuhmacher}\thanksref{t1}\ead[label=e3]{dominic.schuhmacher@stat.unibe.ch}}
\runauthor{L. D\"{u}mbgen, R. Samworth and D. Schuhmacher}
\affiliation{University of Bern, University of Cambridge and University
of Bern}
\address[A]{L. D\"{u}mbgen\\
D. Schuhmacher\\
Institute of Mathematical Statistics\\
\quad and Actuarial Science\\
Alpeneggstrasse 22\\
CH-3012 Bern\\
Switzerland\\
\printead{e1}\\
\phantom{E-mail: }\printead*{e3}}
\address[B]{R. Samworth\\
Statistical Laboratory\\
Centre for Mathematical Sciences\\
Wilberforce Road\\
Cambridge, CB3 0WB\\
United Kingdom\\
\printead{e2}}
\end{aug}

\thankstext{t1}{Supported by the Swiss National Science Foundation.}

\received{\smonth{2} \syear{2010}}
\revised{\smonth{8} \syear{2010}}

%
\begin{abstract}
We study the approximation of arbitrary distributions $P$ on
$d$-di\-mensional space by distributions with log-concave density.
Approximation means minimizing a Kullback--Leibler-type functional. We
show that such an approximation exists if and only if $P$ has finite
first moments and is not supported by some hyperplane. Furthermore we
show that this approximation depends continuously on $P$ with respect
to Mallows distance $D_1(\cdot,\cdot)$. This result implies
consistency of the maximum likelihood estimator of a log-concave
density under fairly general conditions. It also allows us to prove
existence and consistency of estimators in regression models with a
response $Y = \mu(X) + \eps$, where $X$ and $\eps$ are independent,
$\mu(\cdot)$ belongs to a~certain class of regression functions while
$\eps$ is a random error with log-concave density and mean zero.
\end{abstract}

%

\end{frontmatter}

\section{Introduction}
\label{sec:introduction}

Log-concave distributions, that is, distributions with a Lebesgue
density the logarithm of which is concave, are an interesting
nonparametric model comprising many parametric families of
distributions.
\citet{BagnoliBergstrom2005} give an overview of many
interesting properties and applications in econometrics. Indeed, these
distributions have received a lot of attention among statisticians
recently as described in the review by \citet{Walther2009}. The
nonparametric maximum likelihood estimator was studied in the
univariate setting by \citet{Paletal2007},
\citet{Rufibach2006}, \citet{Duembgenetal2007},
\citet{Balabdaouietal2009} and
\citet{DuembgenRufibach2009}. These references contain
characterizations of the estimators, consistency results and explicit
algorithms. Extensions of one or more of these aspects to the
multivariate setting are presented by \citet{Culeetal2010},
\citet{CuleSamworth2010}, \citet{KoenkerMizera2010},
\citet{SereginWellner2011} and \citet
{SchuhmacherDuembgen2010}. Both
\citet{CuleSamworth2010} and \citet{Schuhmacheretal2009}
show that multivariate log-concave
distributions are a very well-behaved nonparametric class. For
instance, moments of arbitrary order are continuous statistical
functionals with respect to weak convergence.

The first aim of the present paper is a deeper understanding of the
approximation scheme underlying the maximum likelihood estimator of a
log-concave density. Let us put this into a somewhat broader context:
let $\hat{Q}_n$ be the empirical distribution of independent random
vectors $X_1$, $X_2, \ldots, X_n$ with distribution $Q$ on a given
open set $\XX\subseteq\R^d$. Suppose that $Q$ has a~density $f$
belonging to a given class $\FF$ of probability densities on $\XX$. The
maximum likelihood estimator of $f$ may be written as
\[
\hat{f}_n = \argmax_{f \in\FF} \int\log(f) \,d\hat{Q}_n
\]
(provided this exists and is unique). Even if $Q$ fails to have a
density within $\FF$, one may view $\hat{f}_n$ as an estimator of the
approximating density
\[
f(\cdot| Q) := \argmax_{f \in\FF} \int\log(f) \,dQ .
\]
In fact, if $Q$ has a density $g \notin\FF$ on $\XX$ such that the
integral $\int g(x) \log g(x) \,dx$ exists in $\R$, one may rewrite
$f(\cdot| Q)$ as the minimizer of the Kullback--Leibler divergence,
\[
D_{\mathrm{KL}}(f,g) = \int\log(g/f)(x) g(x) \,dx,
\]
over all $f \in\FF$. Note the well-known fact that $D_{\mathrm{KL}}(f,g) > 0$
unless $f = g$ almost everywhere. Viewing a maximum likelihood
estimator as an estimator of an approximation within a given model is
common in statistics [see, e.g.,
\citet{Pfanzagl1990}, \citet{Patilea2001},
\citet{Doksumetal2007} and
\citet{CuleSamworth2010}].
\citet{Pfanzagl1990} and \citet{Patilea2001}
show that under suitable regularity conditions on $Q$ and $\FF$, the
estimator $\hat{f}_n$ is consistent with certain large deviation bounds
or rates of convergence, even in the case of misspecified models. To
the best of our knowledge, their results are not directly applicable in
the setting of log-concave densities, which is treated by
\citet{CuleSamworth2010}.
Our ambition is to identify the \textit{largest} possible class of
distributions $Q$ such that $f(\cdot| Q)$ is well defined and
unique. Moreover, we want to show that the mapping $Q \mapsto f(\cdot
| Q)$ is continuous on that class with respect to a \textit{coarse}
topology, ideally the topology of weak convergence.

With these goals in mind, let us tell a short success story about
Grenander's estimator [\citet{Grenander1956}],
also a key example of
\citet{Patilea2001}:
let $\FF_{\mathrm{mon}}$ be the class of all nonincreasing and
left-continuous probability densities on $\XX= (0,\infty)$. Then for
\textit{any} distribution $Q$ on $(0,\infty)$, the maximizer
\[
f_{\mathrm{mon}}(\cdot| Q)
:= \argmax_{f \in\FF_{\mathrm{mon}}} \int_{(0,\infty)} \log f(x) Q(dx)
\]
is well defined and unique. Namely, if $G$ denotes the distribution
function of $Q$, then $f_{\mathrm{mon}}(\cdot| Q)$ is the left-sided
derivative of the smallest concave majorant of $G$ on $(0,\infty)$ [see
\citet{Barlowetal1972}].
With this characterization one can show that for any sequence of
distributions $Q_n$ on $(0,\infty)$ converging weakly to $Q$,
\[
\int_{(0,\infty)} \bigl| f_{\mathrm{mon}}(x | Q_n) - f_{\mathrm{mon}}(x
|
Q) \bigr| \,dx
\to0\qquad (n \to\infty) .
\]
Since\vspace*{1pt} the sequence of empirical measures $\hat{Q}_n$ converges weakly
to $Q$ almost surely, this entails strong consistency of the Grenander
estimator $f(\cdot| \hat{Q}_n)$ in total variation distance.

In the remainder of the present paper we consider the class $\FF$ of
log-con\-cave probability densities on $\XX= \R^d$.
We will show in Section \ref{sec:log-concave approximations} that
$f(\cdot| Q)$ exists and is unique in $L^1(\R^d)$ if and only if
\[
\int\|x\| Q(dx) < \infty
\]
and
\[
Q(H) < 1 \qquad\mbox{for any hyperplane } H \subset\R^d .
\]
Some additional properties of $f(\cdot| Q)$ will be established as
well. We show that the mapping $Q \mapsto f(\cdot| Q)$ is
continuous with respect to
Mallows distance [\citet{Mallows1972}]
$D_1(\cdot,\cdot)$, also known as a Wasserstein, Monge--Kantorovich or
Earth Mover's distance.
Precisely, let $Q$ satisfy the properties just mentioned, and let
$(Q_n)_n$ be a sequence of probability distributions converging to $Q$
in $D_1$; in other words,
%
%
\begin{equation}
\label{eq:Mallows' convergence}
Q_n \to_w Q
\quad\mbox{and}\quad
\int\|x\| Q_n(dx) \to\int\|x\| Q(dx)
\end{equation}
as $n \to\infty$. Then $f(\cdot| Q_n)$ is well defined for
sufficiently large $n$ and
\[
\lim_{n \to\infty} \int\bigl| f(x | Q_n) - f(x | Q)
\bigr|
\,dx
= 0 .
\]
This entails strong consistency of the maximum likelihood estimator
$\hat{f}_n$, because $(\hat{Q}_n)_n$ converges almost surely to $Q$
with respect to Mallows distance $D_1(\cdot,\cdot)$.
In addition we show that $Q \mapsto\max_{f \in\mathcal{F}} \int
\log
(f) \,dQ$ is convex and upper semicontinuous with respect to weak convergence.

In Section \ref{sec:regression problems} we apply these results to the
following type of regression problem: suppose that we observe
independent real random variables $Y_1$, $Y_2, \ldots, Y_n$ such that
\[
Y_i = \mu(x_i) + \eps_i
\]
for given fixed design points $x_1, x_2, \ldots, x_n$ in some set $\XX
$, some unknown regression function $\mu\dvtx\XX\to\R$ and independent,
identically distributed random errors $\eps_i$ with unknown log-concave
density $f$ and mean zero. We will show that a maximum likelihood
estimator of $(\mu,f)$ exists and is consistent under certain
regularity conditions in the following two cases: (i) $\XX= \R^q$ and
$\mu$ is affine (i.e., affine linear); (ii) $\XX= \R$ and $\mu$ is
nondecreasing. These methods are illustrated with a real data set.

Many proofs and technical arguments are deferred to Section \ref
{sec:proofs}. A longer and more detailed version of this paper is the
technical report by \citet{Duembgenetal2010}, referred to as [DSS 2010]
hereafter. It contains all proofs, additional examples and plots, a
detailed description of our algorithms and extensive simulation
studies. There we also indicate potential applications to change-point analyses.

\section{Log-concave approximations}
\label{sec:log-concave approximations}

For a fixed dimension $d \in\mathbb{N}$, let $\Phi= \Phi(d)$ be the
family of concave functions $\phi\dvtx\R^d \to[-\infty,\infty)$ which
are upper semicontinuous and coercive in the sense that
\[
\phi(x) \to- \infty\qquad\mbox{as } \|x\| \to\infty.
\]
In particular, for any $\phi\in\Phi$ there exist constants $a$ and $b
> 0$ such that $\phi(x) \le a - b\|x\|$, so $\int e^{\phi(x)} \,dx$ is
finite. Further let $\QQ= \QQ(d)$ be the family of all probability
distributions $Q$ on $\R^d$. Then we define a log-likelihood-type functional
\[
L(\phi,Q) := \int\phi \,dQ - \int e^{\phi(x)} \,dx + 1
\]
and a profile log-likelihood
\[
L(Q) := \sup_{\phi\in\Phi} L(\phi,Q) .
\]
If, for fixed $Q$, there exists a function $\psi\in\Phi$ such that
$L(\psi,Q) = L(Q) \in\R$, then it will automatically satisfy
\[
\int e^{\psi(x)} \,dx = 1 .
\]
To verify this, note that $\phi+ c \in\Phi$ for any fixed function
$\phi\in\Phi$ and arbitrary $c \in\R$, and
\[
\frac{\partial}{\partial c} L(\phi+ c, Q) = 1 - e^c \int e^{\phi
(x)} \,dx,
\]
if $L(\phi,Q) \in\R$. Thus $L(\phi+c,Q)$ is minimal for $c = - \log
\int e^{\phi(x)} \,dx$.

\subsection{Existence, uniqueness and basic properties}

The next theorem provides a complete characterization of all
distributions $Q \in\QQ$ with real profile log-likelihood $L(Q)$. To
state the result we first define the convex support of a distribution
$Q \in\QQ$ and collect some of its properties.
\begin{Lemma}[(DSS 2010)]
\label{lem:csupp}
For any $Q \in\QQ$, the set
\[
\csupp(Q)
:= \bigcap\{ C \dvtx C \subseteq\R^d \mbox{ closed and
convex}, Q(C) = 1 \}
\]
is itself closed and convex with $Q(\csupp(Q)) = 1$. The following
three proper\-ties of $Q$ are equivalent:

\begin{longlist}[(a)]
\item[(a)] $\csupp(Q)$ has nonempty interior;

\item[(b)] $Q(H) < 1$ for any hyperplane $H \subset\R^d$;

\item[(c)] with $\Leb$ denoting Lebesgue measure on $\R^d$,
\[
\lim_{\delta\downarrow0} \sup\{ Q(C) \dvtx
C \subset\R^d \mbox{ closed and convex}, \Leb(C) \le\delta
\}
< 1 .
\]
\end{longlist}
\end{Lemma}
\begin{Theorem}
\label{thm:existence and uniqueness}
For any $Q \in\QQ$, the value of $L(Q)$ is real if and only if
\[
\int\|x\| Q(dx) < \infty
\quad\mbox{and}\quad
\interior(\csupp(Q)) \ne\varnothing.
\]
In that case, there exists a unique function
\[
\psi= \psi(\cdot| Q) \in\argmax_{\phi\in\Phi} L(\phi
,Q) .
\]
This function $\psi$ satisfies $\int e^{\psi(x)} \,dx = 1$ and
\[
\interior(\csupp(Q)) \subseteq\dom(\psi) := \{x \in\R^d \dvtx
\psi(x)
> - \infty\}
\subseteq\csupp(Q) .
\]
\end{Theorem}
\begin{Remark}[{[Moment (in)equalities]}]
Let $Q \in\QQ$ satisfy the properties stated in Theorem \ref
{thm:existence and uniqueness}. Then the log-density $\psi= \psi
(\cdot
| Q)$ satisfies the following requirements: $L(\psi,Q) = \int\psi
\,dQ \in\R$, and for any function $\Delta\dvtx\R^d \to\R$,
%
%
\begin{equation}
\label{ineq:characterization}
\int\Delta \,dQ \le\int\Delta(x) e^{\psi(x)} \,dx\qquad
\mbox{if } \psi+ t \Delta\in\Phi\qquad\mbox{for some } t > 0 .
\end{equation}
This follows from
\[
\lim_{t \downarrow0} t_{}^{-1} \bigl( L(\psi+ t\Delta,Q) -
L(\psi
,Q) \bigr)
= \int\Delta \,dQ - \int\Delta(x) e^{\psi(x)} \,dx .
\]

Let $P$ be the approximating probability measure with $P(dx) = e^{\psi
(x)} \,dx$. It satisfies the following (in)equalities:
%
%
\begin{eqnarray}
\label{ineq:characterization1}
\int h \,dP & \le& \int h \,dQ\qquad
\mbox{for any convex } h \dvtx\R^d \to(-\infty,\infty] , \\
\label{ineq:characterization2}
\int x P(dx) & = & \int x Q(dx) .
\end{eqnarray}
To verify (\ref{ineq:characterization1}),
let $v \in\R^d$ be a subgradient of $h$ at $0$, that is, $h(x) \ge
h(0) + v^\top x$ for all $x \in\R^d$. Since $\psi(x) \le a - b \|x\|$
for arbitrary $x \in\R^d$ and suitable constants $a$ and $b > 0$, the
function $\Delta:= - h$ satisfies the requirement that $\psi+ t
\Delta\in\Phi$ whenever $0 < t < b/\|v\|$. Hence the asserted
inequality follows from (\ref{ineq:characterization}). The equality for
the first moments follows by setting $h(x) := \pm v^\top x$ for
arbitrary $v \in\R^d$.
\end{Remark}

In what follows let
\begin{eqnarray*}
\QQ^1 = \QQ^1(d) & := & \biggl\{ Q \in\QQ\dvtx\int\|x\| Q(dx) <
\infty\biggr\} , \\
\QQ_o = \QQ_o(d) & := & \{ Q \in\QQ\dvtx\interior(\csupp(Q))
\ne
\varnothing\} .
\end{eqnarray*}
Thus $L(Q) \in\R$ if and only if $Q \in\QQ_o \cap\QQ^1$. Moreover,
the proof of Theorem~\ref{thm:existence and uniqueness} shows that
\[
L(Q) = \cases{
- \infty, &\quad for $Q \in\QQ\setminus\QQ^1 $, \cr
+ \infty, &\quad for $Q \in\QQ^1 \setminus\QQ_o $.}
\]
\begin{Remark}[(Affine equivariance)]
\label{rem:equivariance}
Suppose that $Q \in\QQ_o \cap\QQ^1$. For arbitrary vectors $a \in
\R
^d$ and nonsingular, real $d \times d$ matrices $B$ define $Q_{a,B}$ to
be the distribution of $a + BX$ when $X$ has distribution $Q$. Then
$Q_{a,B} \in\QQ_o\cap\QQ^1$, too, and elementary considerations
reveal that
\[
L(Q_{a,B}) = L(Q) - {\log}|{\det B}|
\]
and
\[
\psi(x | Q_{a,B}) = \psi\bigl( B^{-1}(x - a) | Q\bigr) -
{\log}|{\det B}|\qquad
\mbox{for } x \in\R^d .
\]
\end{Remark}
\begin{Remark}[(Convexity, DSS 2010)]
\label{rem:convexity}
The profile log-likelihood $L$ is convex on $\QQ^1$. Precisely, for
arbitrary $Q_0, Q_1 \in\QQ^1$ and $0 < t < 1$,
\[
L \bigl( (1 - t)Q_0 + tQ_1 \bigr) \le(1 - t) L(Q_0) + t L(Q_1) .
\]
The two sides are equal and real if and only if $Q_0, Q_1 \in\QQ_o
\cap\QQ^1$ with $\psi(\cdot| Q_0) = \psi(\cdot| Q_1)$.
\end{Remark}
\begin{Remark}[(Concave majorants, DSS 2010)]
\label{rem:concave majorants}
Let $\psi= \psi(\cdot| Q)$ for a distribution $Q \in\QQ_o \cap
\QQ^1$. For any open set $U \subset\R^d$ there exists a (pointwise)
minimal function $\psi_U \in\Phi$ such that $\psi_U \ge\psi$ on
$\R^d
\setminus U$. In particular, $\psi_U \le\psi$ with equality on $\R^d
\setminus U$. One can also show that $\psi_U$ is the pointwise infimum
of all affine functions $\phi$ such that $\phi\ge\psi$ on $\R^d
\setminus U$. If $\operatorname{supp}(Q)$ denotes the smallest closed set $A
\subseteq\R^d$ with $Q(A) = 1$, then
\[
\psi= \psi_{\R^d \setminus\mathrm{supp}(Q)} .
\]
Furthermore, suppose that $Q$ has a density $g$ on an open set $U$ such
that $\psi> \log g$ on this set. Then
\[
\psi= \psi_U .
\]
\end{Remark}

\subsection{The one-dimensional case}
\label{subsec:dimension one}

For the case of $d = 1$ one can generalize Theorem 2.4 of
\citet{DuembgenRufibach2009}
as follows: for a function $\phi\in\Phi(1)$ let
\[
\SSS(\phi) := \bigl\{ x \in\dom(\phi) \dvtx
\phi(x) > 2^{-1} \bigl( \phi(x - \delta) + \phi(x + \delta) \bigr)
\mbox{ for all } \delta> 0 \bigr\}.
\]
The log-concave approximation of a distribution on $\R$ can be
characterized in terms of distribution functions only:
\begin{Theorem}
\label{thm:characterization for d=1}
Let $Q$ be a nondegenerate distribution on $\R$ with finite first
moment and distribution function $G$. Let $F$ be a distribution
function with log-density $\phi\in\Phi$. Then $\phi= \psi(\cdot
| Q)$ if and only if
\[
\int_{-\infty}^\infty\bigl( F(t) - G(t) \bigr) \,dt
= 0
\]
and
\[
\int_{-\infty}^x \bigl( F(t) - G(t) \bigr) \,dt
\cases{
\le0, &\quad for all $x \in\R$, \cr
= 0, &\quad for all $x \in\SSS(\phi) $.}
\]
\end{Theorem}
\begin{Remark}[(DSS 2010)]
One consequence of this theorem is that the c.d.f. $F$ of $\psi(\cdot
| Q)$ follows the c.d.f. $G$ of $Q$ quite closely in that
\[
G(x -) \le F(x) \le G(x)\qquad
\mbox{for arbitrary } x \in\SSS(\psi(\cdot| Q)) .
\]
\end{Remark}
\begin{Example}
\label{ex:ProjectionLC1}
Let $Q$ be a rescaled version of Student's distribution $t_2$ with
density and distribution function
\[
g(x) = 2^{-1} (1 + x^2)^{-3/2}
\quad\mbox{and}\quad
G(x) = 2^{-1}\bigl(1 + (1 + x^2)^{-1/2} x\bigr) ,
\]
respectively. The best approximating log-concave distribution is the
Laplace distribution with density and distribution function
\[
f(x) = 2^{-1} e^{-|x|}
\quad\mbox{and}\quad
F(x) = \cases{
f(x), &\quad for $x \le0 $, \cr
1 - f(x), &\quad for $x \ge0 $,}
\]
respectively. To verify this claim, note that by symmetry it suffices
to show that
\[
\int_{-\infty}^x \bigl( F(t) - G(t) \bigr) \,dt
\cases{
\le0, &\quad for $x \le0$, \cr
= 0, &\quad for $x = 0 $.}
\]
Indeed the integral on the left-hand side equals
\[
2^{-1} \bigl( \exp(x) - x - (1 + x^2)^{1/2} \bigr)
\]
for all $x \le0$. Clearly this expression is zero for $x = 0$, and
elementary considerations show that it is nonpositive for all $x \le
0$. Numerical calculations reveal that $|F - G|$ is smaller than $0.04$
everywhere.
\end{Example}
\begin{Remark}
Let $Q \in\QQ_o \cap\QQ^1$ such that $Q(a,b) = 0$ for some bounded
interval $(a,b) \subset\csupp(Q)$. Then $\psi= \psi(\cdot| Q)$
is linear on $[a,b]$. This follows from Remark \ref{rem:concave
majorants}, applied to $U = (a,b)$. Note that $\psi(a) > - \infty$ and
$\psi(b) > - \infty$, because otherwise $\psi\equiv- \infty$ on
$(-\infty,a]$ or on $[b,\infty)$. But this would be incompatible with
$\int\psi \,dQ \in\R$, because both $Q((-\infty,a])$ and
$Q([b,\infty
))$ are positive.
\end{Remark}
\begin{Remark}[(DSS 2010)]
\label{rem:examples for d=1}
Suppose that $Q$ has a continuous but not log-concave density $g$.
Nevertheless one can say the following about the approximating
log-density $\psi= \psi(\cdot| Q)$:

\begin{longlist}
\item Suppose that $\log g$ is concave on an interval $(-\infty
,a]$ with $g(a) > 0$ and $\psi(a) \le\log g(a)$. Then there exists a
point $a' \in[-\infty,a]$ such that $\psi$ is linear on $(a',a]$ and
$\psi= \log g$ on $(-\infty,a']$.
\item Suppose that $\log g$ is differentiable everywhere, convex
on a bounded interval $[a,b]$ and concave on both $(-\infty,a]$ and
$[b,\infty)$. Then there exist points $a' \in(-\infty,a]$ and $b'
\in
[b,\infty)$ such that $\psi$ is linear on $[a',b']$ while $\psi=
\log
g$ on $(-\infty,a'] \cup[b',\infty)$.
\item Suppose that $\log g$ is convex on an interval $(-\infty
,a]$ such that $- \infty< \log g(a) \le\psi(a)$. Then $\psi$ is
linear on $(-\infty,a]$.
\end{longlist}
\end{Remark}
\begin{Example}
Let us illustrate part (ii) of Remark \ref{rem:examples for d=1} with a
numeri\-cal example. Figure \ref{fig:ProjectionLC2} shows the bimodal
%
%
\begin{figure}

\includegraphics{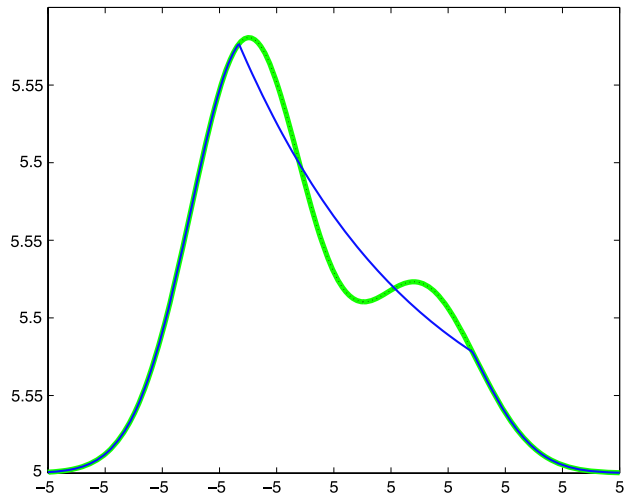}

\caption{Density of a Gaussian mixture and its log-concave approximation.}
\label{fig:ProjectionLC2}
\end{figure}
density $g$
(green/dotted line)
of the Gaussian mixture $Q = 0.7 \cdot\mathcal{N}(-1.5,1) + 0.3 \cdot
\mathcal{N}(1.5,1)$ together with its log-concave approximation $f =
f(\cdot| Q)$ (blue line).
As predicted, there exists an interval $[a',b']$ such that $f = g$ on
$\R\setminus(a',b')$ and $\log f$ is linear on $[a',b']$.
\end{Example}

\subsection{Continuity in $Q$}

For the applications to regression problems to follow we need to
understand the properties of both $Q \mapsto L(Q)$ and $Q \mapsto\psi
(\cdot| Q)$ on $\QQ^1 \cap\QQ_o$. Our first hope was that both
mappings would be continuous~with respect to the weak topology. It
turned out, however, that we need a somewhat stronger notion of
convergence, namely, convergence with respect to Mallows distance
$D_1$ which is defined as follows: for two probability distributions
$Q, Q' \in\QQ^1$,
\[
D_1(Q,Q') := \inf_{(X,X')} \Ex\|X - X'\| ,
\]
where the infimum is taken over all pairs $(X,X')$ of random vectors $X
\sim Q$ and $X' \sim Q'$ on a common probability space. It is well
known that the infimum in $D_1(Q,Q')$ is a minimum. The distance $D_1$
is also known as Wasserstein, Monge--Kantorovich or Earth Mover's
distance. An alternative representation due to
\citet{KantorovichRubinstein1958}
is
\[
D_1(Q,Q')
= \sup_{h \in\HH_L}
\biggl| \int h \,dQ - \int h \,dQ' \biggr|,
\]
where $\HH_L$ consists of all $h\dvtx\R^d \to\R$ such that $|h(x) - h(y)|
\le\|x-y\|$ for all $x,y \in\R^d$. Moreover, for a sequence $(Q_n)_n$
in $\QQ^1$, it is known that (\ref{eq:Mallows' convergence}) is
equivalent to $D_1(Q_n,Q) \to0$ as $n \to\infty$
[\citet{Mallows1972}, \citet{BickelFreedman1981}]. In case
of $d = 1$, the optimal coupling of $Q$
and $Q'$ is given by the quantile transformation: if $G$ and $G'$
denote the respective distribution functions, then
\[
D_1(Q,Q') = \int_0^1 | G^{-1}(u) - {G'}^{-1}(u) | \,du
= \int_{-\infty}^{\infty} | G(x) - G'(x) | \,dx .
\]
A good starting point for more detailed information on Mallows
distance is Chapter 7 of
\citet{Villani2003}.

Before presenting the main results of this section we mention two
useful facts about the convex support of distributions.
\begin{Lemma}
\label{lem:weak semicontinuity of csupp}
Given a distribution $Q \in\QQ$, a point $x \in\R^d$ is an interior
point of $\csupp(Q)$ if and only if
\[
h(Q,x) := \sup\{ Q(C) \dvtx C \subset\R^d \mbox{ closed and convex},
x \notin\interior(C) \}
< 1 .
\]
Moreover, if $(Q_n)_n$ is a sequence in $\QQ$ converging weakly to
$Q$, then
\[
\limsup_{n \to\infty} h(Q_n,x) \le h(Q,x)
\qquad\mbox{for any } x \in\R^d .
\]
\end{Lemma}

This lemma implies that the set $\QQ_o$ is an open subset of $\QQ$ with
respect to the topology of weak convergence. The supremum $h(Q,x)$ is a
maximum over closed halfspaces and is related to Tukey's halfspace
depth [\citet{DonohoGasko1992}, Section 6]. For a proof of
Lemma \ref{lem:weak semicontinuity of csupp} we refer to [DSS 2010]. Now
we are ready to state the main results of this section.
\begin{Theorem}[(Weak upper semicontinuity)]
\label{thm:weak upper semicontinuity}
Let $(Q_n)_n$ be a sequence of distributions in $\QQ_o$ converging
weakly to some $Q \in\QQ_o$. Then
\[
\limsup_{n \to\infty} L(Q_n) \le L(Q) .
\]
Moreover, $\liminf_{n \to\infty} L(Q_n) < L(Q)$ if and only if
\[
\limsup_{n \to\infty} \int\|x\| Q_n(dx) > \int\|x\|
Q(dx) .
\]
\end{Theorem}

This result already entails continuity of $L(\cdot)$ on $\QQ_o\cap
\QQ
^1$ with respect to Mallows distance $D_1$. The next theorem extends
this result to $L \dvtx\QQ^1 \to(-\infty,\infty]$:
\begin{Theorem}[(Continuity with respect to Mallows distance $D_1$)]
\label{thm:mallows-continuity}
Let $(Q_n)_n$ be a sequence of distributions in $\QQ^1$ such that
$\lim
_{n \to\infty} D_1(Q_n,Q) = 0$ for some $Q \in\QQ^1$. Then
\[
\lim_{n \to\infty} L(Q_n) = L(Q) .
\]
In case of $Q \in\QQ_o \cap\QQ^1$, the probability densities $f :=
\exp\mbox{${}\circ{}$} \psi(\cdot| Q)$ and $f_n := \exp\mbox{${}\circ{}$} \psi
(\cdot
| Q_n)$ are well defined for sufficiently large $n$ and satisfy
\begin{eqnarray*}
\lim_{n \to\infty, x \to y} f_n(x) & = & f(y)
\qquad\mbox{for all } y \in\R^d \setminus\partial\{f > 0\} , \\
\limsup_{n \to\infty, x \to y} f_n(x) & \le& f(y)
\qquad\mbox{for all } y \in\partial\{f > 0\}, \\
\lim_{n \to\infty} \int| f_n(x) - f(x) | \,dx
& = & 0 .
\end{eqnarray*}
\end{Theorem}
\begin{Remark}[(Stronger modes of convergence)]
The convergence of $(f_n)_n$ to $f$ in total variation distance may be
strengthened considerably. It follows from recent results of
\citet{CuleSamworth2010} or
\citet{Schuhmacheretal2009}
that $(f_n)_n \to f$ uniformly on arbitrary closed subsets of $\R^d
\setminus\operatorname{disc}(f)$, where $\operatorname{disc}(f)$ is the set of
discontinuity points of $f$. The latter set is contained in the
boundary of the convex set $\{f > 0\}$, hence a null set with respect
to Lebesgue measure. Moreover, there exists a~number $\eps(f) > 0$
such that
\[
\lim_{n \to\infty}
\int e^{\eps(f) \|x\|} | f_n(x) - f(x) | \,dx
= 0 .
\]
More generally,
\[
\lim_{n \to\infty}
\int e^{A(x)} | f_n(x) - f(x) | \,dx
= 0
\]
for any sublinear function $A \dvtx\R^d \to\R$ such that $\lim_{\|x\|
\to
\infty} e^{A(x)} f(x) = 0$.
\end{Remark}

\section{Applications to regression problems}
\label{sec:regression problems}

Now we consider the regression setting described in the
\hyperref[sec:introduction]{Introduction}
with observations $Y_i = \mu(x_i) + \eps_i$, $1 \le i \le n$, where the
$x_i \in\XX$ are given fixed design points, $\mu\dvtx\XX\to\R$ is an
unknown regression function, and the $\eps_i$ are independent random
errors with mean zero and unknown distribution $Q$ on $\R$ such that
$\psi= \psi(\cdot| Q)$ is well defined. The regression function
$\mu$ is assumed to belong to a given family $\MM$ with the property that
\[
m + c \in\MM\qquad\mbox{for arbitrary } m \in\MM, c \in\R.
\]

\subsection{Maximum likelihood estimation}
\label{ssec:maximum likelihood estimation}

We propose to estimate $(\psi,\mu)$ by a maximizer of
\[
\hat{\Lambda}(\phi,m) := \frac{1}{n} \sum_{i=1}^n \phi\bigl(Y_i -
m(x_i)\bigr) - \int e^{\phi(x)} \,dx + 1
\]
over all $(\phi,m) \in\Phi\times\MM$. Note that $\hat{\Lambda
}(\phi
,m)$ remains unchanged if we replace $(\phi,m)$ with $(\phi(\cdot+ c),
m+c)$ for an arbitrary $c \in\R$. For fixed $m$, the maximizer $\hat
{\phi} = \hat{\phi}_m$ of $\hat{\Lambda}(\cdot,m)$ over $\Phi$ will
automatically satisfy $\int\exp(\hat{\phi}(x)) \,dx = 1$ and $\int x
\exp(\hat{\phi}(x)) \,dx = n^{-1} \sum_{i=1}^n (Y_i - m(x_i))$. Thus
if $(\hat{\phi},\hat{m})$ maximizes $\hat{\Lambda}(\cdot,\cdot)$ over
$\Phi\times\MM$, then
\[
(\hat{\psi},\hat{\mu}) := \bigl(\hat{\phi}(\cdot+ c), \hat{m} + c\bigr)
\qquad\mbox{with }
c := \frac{1}{n}\sum_{i=1}^n \bigl(Y_i - \hat{m}(x_i)\bigr)
\]
maximizes $\hat{\Lambda}(\phi,m)$ over all $(\phi,m) \in\Phi
\times\MM
$ satisfying the additional constraint that $\exp\mbox{${}\circ{}$} \phi$ defines
a probability density with mean zero.

Define $\bs{x} := (x_i)_{i=1}^n$ and $m(\bs{x}) := (m(x_i))_{i=1}^n$.
Then we may write
\[
\hat{\Lambda}(\phi,m) = L\bigl(\phi,\hat{Q}_{m(\bs{x})}\bigr)
\]
with the empirical distributions
\[
\hat{Q}_{\bs{v}} := \frac{1}{n} \sum_{i=1}^n \delta_{Y_i - v_i}
\]
for $\bs{v} = (v_i)_{i=1}^n \in\R^n$. Thus our procedure aims to find
\[
(\hat{\phi}, \hat{m})
\in\argmax_{(\phi,m) \in\Phi\times\MM} L\bigl(\phi,\hat
{Q}_{m(\bs
{x})}\bigr) ,
\]
and this representation is our key to proving the existence of $(\hat
{\psi},\hat{\mu})$. Before doing so we state a simple inequality of
independent interest, which follows from Jensen's inequality and
elementary considerations:
\begin{Lemma}[(DSS 2010)]
\label{lem:life is simple}
For any distribution $Q \in\QQ^1(1)$,
\[
L(Q) \le- \log\biggl( 2 \int|x - \operatorname{Med}(Q)| Q(dx)
\biggr)
\le- \log\biggl( \int|x - \mu(Q)| Q(dx) \biggr) ,
\]
where $\operatorname{Med}(Q)$ is a median of $Q$ while $\mu(Q)$ denotes its
mean $\int x Q(dx)$.
\end{Lemma}
\begin{Theorem}[(Existence in regression)]
\label{thm:existence regression}
Suppose that the set $\MM(\bs{x}) := \{ m(\bs{x}) \dvtx m \in\MM
\} \subset\R^n$ is closed and does not contain $\bs{Y} :=
(Y_i)_{i=1}^n$. Then there exists a maximizer $(\hat{\phi},\hat{m})$ of
$\hat{\Lambda}(\phi,m)$ over all $(\phi,m) \in\Phi\times\MM$.
\end{Theorem}

The constraint $\bs{Y} \notin\MM(\bs{x})$ excludes situations with
perfect fit. In that case, the Dirac measure $\delta_0$ would be the
most plausible error distribution.
\begin{Example}[(Linear regression)]
Let $\XX= \R^q$, and let $\MM$ consist of all affine functions on
$\R
^q$. Then $\MM(\bs{x})$ is the column space of the design matrix
\[
\bs{X} = \left[\matrix{
1 & 1 & \cdots& 1 \cr
x_1 & x_2 & \cdots& x_n}\right]^\top
\in\R^{n \times(q+1)} ,
\]
hence a linear subspace of $\R^n$. Consequently there exists a
maximizer $(\hat{\phi},\hat{m})$ of $\hat{\Lambda}$ over $\Phi
\times
\MM$, unless $\bs{Y} \in\MM(\bs{x})$.
\end{Example}
\begin{Example}[(Isotonic regression)]
Let $\XX$ be some interval on the real line, and let $\MM$ consist of
all isotonic functions $m \dvtx\XX\to\R$. Then the set $\MM(\bs{x})$ is
a closed convex cone in $\R^n$. Here the condition that $\bs{Y} \notin
\MM(\bs{x})$ is equivalent to the existence of two indices $i,j \in\{
1,2,\ldots,n\}$ such that $x_i \le x_j$ but $Y_i > Y_j$.
\end{Example}

\subsubsection*{Fisher consistency}
The maximum likelihood estimator $(\hat{\psi},\hat{\mu})$ need not be
unique in general. Nevertheless we will prove it to be consistent under
certain regularity conditions. A key point here is Fisher consistency
in the following sense: note that the expectation measure of the
empirical distribution $\hat{Q}_{m(\bs{x})}$ equals
\[
\Ex\hat{Q}_{m(\bs{x})}
= \frac{1}{n} \sum_{i=1}^n Q \star\delta_{\mu(x_i) - m(x_i)}
= Q \star R_{(\mu- m)(\bs{x})}
\]
with
\[
R_{\bs{v}} := \frac{1}{n} \sum_{i=1}^n \delta_{v_i} .
\]
But
\[
L\bigl(Q\star R_{(\mu- m)(\bs{x})}\bigr) \le L(Q)
\]
with equality if and only if $\mu- m$ is constant on $\{
x_1,x_2,\ldots,x_n\}$. This follows from a more general inequality
which is somewhat reminiscent of
Anderson's lemma [\citet{Anderson1955}]:
\begin{Theorem}
\label{thm:convolution}
Let $Q \in\QQ_o(d) \cap\QQ^1(d)$ and $R \in\QQ^1(d)$. Then
$Q\star R
\in\QQ_o \cap\QQ^1$ and
\[
L(Q\star R) \le L(Q) .
\]
Equality holds if and only if $R = \delta_a$ for some $a \in\R^d$.
\end{Theorem}

\subsection{Consistency}

In this subsection we consider a triangular scheme of independent
observations $(x_{ni},Y_{ni})$, $1 \le i \le n$, with fixed design
points $x_{ni} \in\XX_n$ and
\[
Y_{ni} = \mu_n(x_{ni}) + \eps_{ni} ,
\]
where $\mu_n$ is an unknown regression function in $\MM_n$ and $\eps
_{n1}, \eps_{n2}, \ldots, \eps_{nn}$ are unobserved independent random
errors with mean zero and unknown distribution $Q_n \in\QQ_o(1) \cap
\QQ^1(1)$. Two basic assumptions are:

\begin{longlist}[(A.2)]
\item[(A.1)] $\MM_n(\bs{x}_n)$ is a closed subset of $\R^n$ for every
$n \in\mathbb{N}$;
\item[(A.2)] $D_1(Q_n,Q) \to0$ for some distribution $Q \in\QQ_o(1)
\cap\QQ^1(1)$.
\end{longlist}

We write $(\hat{\psi}_n, \hat{\mu}_n)$ for a maximizer of $L(\phi
,\hat
{Q}_{n,m})$ over all pairs $(\phi,m) \in\Phi\times\MM_n$ such that
$\int e^{\phi(x)} \,dx = 1$ and $\int x e^{\phi(x)} \,dx = 0$, where
$\hat{Q}_{n,m}$ stands for the empirical distribution of the residuals
$Y_{ni} - m(x_{ni})$, $1 \le i \le n$. We also need to consider its
expectation measure
\[
Q_{n,m} := \Ex\hat{Q}_{n,m} = Q_n \star R_{(\mu_n - m)(\bs
{x}_n)} .
\]
Furthermore we write
\[
\|\bs{v}\|_n := \frac{1}{n} \sum_{i=1}^n |v_i|\qquad
\mbox{for } \bs{v} = (v_i)_{i=1}^n \in\R^n .
\]
It is also convenient to metrize weak convergence. In Theorem \ref
{thm:consistency regression} below we utilize the bounded Lipschitz
distance: for probability distributions $Q,Q'$ on the real line let
\[
D_{\mathrm{BL}}(Q,Q') := \sup_{h \in\HH_{\mathrm{BL}}} \biggl| \int h \, d(Q - Q')
\biggr| ,
\]
where $\HH_{\mathrm{BL}}$ is the family of all functions $h \dvtx\R\to[-1,1]$
such that $| h(x) - h(y) | \le|x - y|$ for all $x,y \in\R$.
\begin{Theorem}[(Consistency in regression)]
\label{thm:consistency regression}
Let assumptions \textup{(A.1)} and \textup{(A.2)} be satisfied. Suppose further
that:
\begin{longlist}[(A.2)]
\item[(A.3)] for arbitrary fixed $c > 0$,
\[
\sup_{m \in\MM_n \dvtx \|(m - \mu_n)(\bs{x}_n)\|_n \le c}
D_{\mathrm{BL}}(\hat{Q}_{n,m},Q_{n,m})
\to_p 0 .
\]
Then, with $f_n := \exp\mbox{${}\circ{}$} \psi(\cdot| Q_n)$ and $\hat{f}_n
:= \exp\mbox{${}\circ{}$} \hat{\psi}_n$, the maximum likelihood estimator
$(\hat
{f}_n,\hat{\mu}_n)$ of $(f_n,\mu_n)$ exists with asymptotic probability
one and satisfies
\[
\int| \hat{f}_n(x) - f_n(x) | \,dx \to_p 0 ,\qquad
\| (\hat{\mu}_n - \mu_n)(\bs{x}_n) \|_n \to_p 0 .
\]
\end{longlist}
\end{Theorem}

We know already that assumption (A.1) is satisfied for multiple linear
regression and isotonic regression. Assumption (A.2) is a
generalization of assuming a fixed error distribution for all sample
sizes. The crucial point, of course, is assumption (A.3). In our two
examples it is satisfied under mild conditions:
\begin{Theorem}[(Linear regression)]
\label{thm:consistency linear regression} Let $\MM_n$ be the family of
all affine functions on $\XX_n := \R ^{q(n)}$. If assumption
\textup{(A.2)} is satisfied, then \textup{(A.3)} follows from
\[
\lim_{n\to\infty} q(n)/n = 0 .
\]
\end{Theorem}
\begin{Theorem}[(Isotonic regression)]
\label{thm:consistency isotonic regression}
Let $\MM_n$ be the set of all nondecreasing functions on an interval
$\XX_n \subseteq\R$. If assumption \textup{(A.2)} holds true, then
\textup{(A.3)} follows from
\[
\| \mu_n(\bs{x}_n) \|_n = O(1) .
\]
\end{Theorem}

The proof of Theorem \ref{thm:consistency linear regression} is given
in Section \ref{sec:proofs}. For the proof of Theorem~\ref{thm:consistency isotonic regression},
which uses similar ideas and an
additional approximation argument, we refer to [DSS 2010].

\subsection{Algorithms and numerical results}

Computing the maximum likelihood estimator $(\hat{\psi}, \hat{\mu})$
from Section \ref{ssec:maximum likelihood estimation} turns out to
be a rather difficult task, because the function $\hat{\Lambda}$ can
have multiple local maxima. In [DSS 2010] we discuss strengths and
weaknesses of three different algorithms, including an alternating and
a stochastic search algorithm. The third procedure, which is highly
successful in the case of linear regression,\vspace*{1pt} is global maximization of
the profile log-likelihood $\hat{\Lambda}(\theta) := \max_{\phi\in
\Phi
} \hat{\Lambda}(\phi,m_{\theta})$, where $m_{\theta}(x) = \theta
^{\top}
x$ for every $x \in\R^q$, by means of differential
evolution [\citet{Priceetal2005}].

Extensive simulation studies in [DSS 2010] suggest that $(\hat{\psi},
\hat
{\mu})$ provides rather accurate estimates even if $n$ is only
moderately large. For various skewed error distributions, $\hat{\mu}$
may be considerably better than the corresponding least squares
estimator. As an example consider the simple linear regression model
with observations
\[
Y_i = c + \theta X_i + \eps_i,\qquad 1 \le i \le n := 100,
\]
where $X_1, \ldots, X_n$ are independent design points from the
$\operatorname{Unif}[0,3]$ distribution and $\eps_1, \ldots, \eps_n$ are independent
errors from a centered gamma distribution with shape parameter $r$ and
variance $1$. Note that the distribution of $(\hat{\psi},\hat{\theta
}-\theta)$ does not depend on $c$ or $\theta$. Monte Carlo estimation
of the root mean squared error based on 1000 simulations of this model
gives 0.023 for the estimator $\hat{\theta}$ versus 0.118 for the least
squares estimator of $\theta$ if $r=1$, and 0.095 versus 0.113 for the
same comparison if $r=3$.

\subsection{A data example}

A familiar task in econometrics is to model expenditure ($Y$) of
households as a function of their income ($X$). Not only the mean curve
(Engel curve) but also quantile curves play an important role. A
related application are growth charts in which, for instance, $X$ is
the age of a newborn or infant and $Y$ is its height or weight.

\begin{figure}

\includegraphics{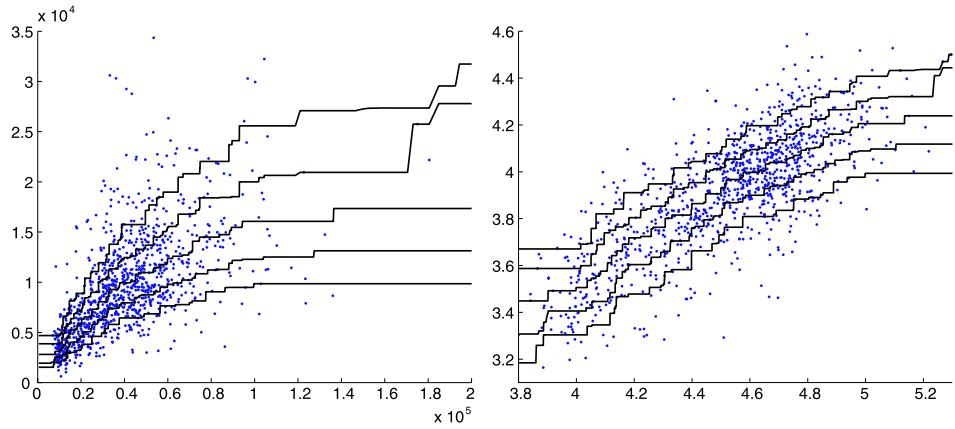}

\caption{UK household data, raw (left) and log-transformed
(right),
with isotonic quantile curves.}
\label{fig:UK_data}
\end{figure}

We applied our methods to a survey of $n = 7125$ households in the
United Kingdom in 1973 (data courtesy of W. H\"ardle, HU Berlin). The
two variables we considered were annual income ($X_{\mathrm{raw}}$) and
annual expenditure for food ($Y_{\mathrm{raw}}$). Figure \ref{fig:UK_data}
shows scatter plots of the raw and log-transformed data. To enhance
visibility we only show a random subsample of size $n' = 1000$. In
addition, isotonic quantile curves $x \mapsto\hat{q}_\beta(x)$ are
added for $\beta= 0.1$, $0.25$, $0.5$, $0.75$, $0.9$ (based on all
observations). These pictures show clearly that the raw data are
heteroscedastic, whereas for the log-transformed data, $(X_i,Y_i) =
(\log_{10} X_{{\mathrm{raw}},i}, \log_{10} Y_{{\mathrm{raw}},i})$, an additive
model seems appropriate.

%

Interestingly, neither linear nor quadratic nor cubic regression yield
convincing fits to these data. Polynomial regression of degree four or
cubic splines with knot points at, say, $4.1$, $4.3$, $4.5$, $4.7$,
$4.9$ seem to fit the data quite well. Moreover, exact Monte Carlo
goodness-of-fits test, assuming the regression function to be a cubic
spline and based on a Kolmogorov--Smirnov statistic applied to
studentized residuals, revealed the regression errors $\eps_i$ to be
definitely non-Gaussian.

%
\begin{figure}[b]

\includegraphics{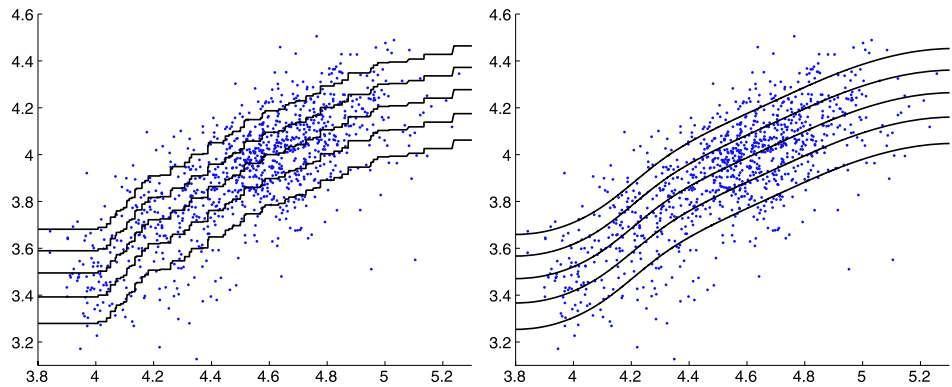}

\caption{Log-transformed UK household data with isotonic fits
(left)
and spline fits (right) from our additive model.}
\label{fig:UK_fits}
\end{figure}

Figure \ref{fig:UK_fits} shows the data and estimated $\beta$-quantile
curves for $\beta= 0.1$, $0.25$, $0.5$, $0.75$, $0.9$, based on our
additive regression model. Note that the estimated $\beta$-quantile
curve is simply the estimated mean curve plus the $\beta$-quantile of
the estimated error distribution. On the left-hand side, we only
assumed $\mu$ to be nondecreasing, on the right-hand side we fitted
the aforementioned spline model. In both cases the fitted quantile
curves are similar to the quantile curves in Figure~\ref{fig:UK_data}
but with fewer irregularities such as big jumps which may be artifacts
due to sampling error.

\section{Proofs}
\label{sec:proofs}

For the proof of Theorem \ref{thm:existence and uniqueness} we need an
elementary bound for the Lebesgue measure of level sets of log-concave
distributions:
\begin{Lemma}[(DSS 2010)]
\label{lem:level sets}
Let $\phi\in\Phi$ be such that $\int e^{\phi(x)} \,dx = 1$. For real
$t$ define the level set $D_t := \{x \in\R^d \dvtx\phi(x) \ge t\}$. Then
for $r < M \le\max_{x \in\R^d} \phi(x)$,
\[
\Leb(D_r) \le(M - r)^d e^{-M} \Big/ \int_0^{M - r} t^d e^{-t}
\,dt .
\]
\end{Lemma}

Another key ingredient for the proofs of Theorems \ref{thm:existence
and uniqueness} and \ref{thm:mallows-continuity} is a lemma on
pointwise limits of sequences in $\Phi$:
\begin{Lemma}[(DSS 2010)]
\label{lem:limits of sequences in GG}
Let $\bar{\phi}$ and $\phi_1, \phi_2, \phi_3, \ldots$ be
functions in
$\Phi$ such that $\phi_n \le\bar{\phi}$ for all $n \in\mathbb{N}$.
Further suppose that the set
\[
C := \Bigl\{ x \in\R^d \dvtx\liminf_{n \to\infty} \phi_n(x)
> -
\infty\Bigr\}
\]
is nonempty. Then there exist a subsequence $(\phi_{n(k)})_k$ of
$(\phi
_n)_n$ and a function $\phi\in\Phi$ such that $C \subseteq\dom
(\phi)
= \{\phi> - \infty\}$ and
\begin{eqnarray*}
\lim_{k \to\infty, x \to y} \phi_{n(k)}(x) & = & \phi(y)\qquad
\mbox{for all } y \in\interior(\dom(\phi)) , \\
\limsup_{k \to\infty, x \to y} \phi_{n(k)}(x) & \le& \phi(y)
\le\bar{\phi}(y)\qquad
\mbox{for all } y \in\R^d .
\end{eqnarray*}
\end{Lemma}
\begin{pf*}{Proof of Theorem \ref{thm:existence and uniqueness}}
Suppose first that $\int\|x\| Q(dx) = \infty$. Since any $\phi\in
\Phi$ is majorized by $x \mapsto a - b \|x\|$ for suitable constants
$a$ and $b > 0$, this entails that $L(Q) = - \infty$.

Second, suppose that $\int\|x\| Q(dx) < \infty$ but $\interior
(\csupp(Q)) = \varnothing$. According to Lemma \ref{lem:csupp}, the
latter fact is equivalent to $Q(H) = 1$ for some hyperplane $H \subset
\R^d$. For $c \in\R$ define a function $\phi_c \in\Phi$ via $\phi
_c(x) := c - \|x\|$ for $x \in H$ and $\phi_c(x) := - \infty$ for $x
\notin H$. Then $L(\phi_c,Q) = c - \int\|x\| Q(dx) + 1 \to
\infty$
as $c \to\infty$.

For the remainder of this proof suppose that $\int\|x\| Q(dx) < \infty
$ and that $\csupp(Q)$ has nonempty interior. Since the concave
function $h(x) = -\|x\|$ satisfies $\int h \,dQ > -\infty$, we have
$L(Q) > -\infty$. When maximizing $L(\phi,Q)$ over all $\phi\in\Phi$
we may and do restrict our attention to functions $\phi\in\Phi$ such
that $\int e^{\phi(x)} \,dx = 1$ (see end of Section \ref
{sec:introduction}) and $\dom(\phi) = \{\phi> - \infty\} \subseteq
\csupp(Q)$. For if $\dom(\phi) \not\subset\csupp(Q)$, replacing
$\phi
(x)$ with $-\infty$ for all $x \notin\csupp(Q)$ would also increase
$L(\phi,Q)$ strictly. Let $\Phi(Q)$ be the family of all $\phi\in
\Phi
$ with these properties.

Now we show that $L(Q) < \infty$. Suppose that $\phi\in\Phi(Q)$ is
such that $M := \max_{x \in\mathbb{R}^d} \phi(x) > 0$. With $D_t :=
\{
\phi\ge t\}$ and for $c > 0$ we get the bound
\begin{eqnarray*}
L(\phi,Q) &=& \int\phi \,dQ
\le -cM Q(\R^d \setminus D_{-cM}) + M Q(D_{-cM}) \\
& = & - (c+1) M \biggl( \frac{c}{c+1} - Q(D_{-cM}) \biggr) .
\end{eqnarray*}
According to Lemma \ref{lem:level sets},
\begin{eqnarray*}
\Leb(D_{-cM})
& \le& (1+c)^d M^d e^{-M} \Big/ \int_0^{(1+c)M} t^d e^{-t} \,dt \\
& = & (1+c)^d M^d e^{-M} / \bigl(d! + o(1)\bigr)
\to0
\end{eqnarray*}
as $M \to\infty$ for any fixed $c > 0$. But Lemma \ref{lem:csupp}
entails that for sufficiently large $c$ and sufficiently small $\delta
> 0$,
\[
\sup\{ Q(C) \dvtx
C \subset\R^d \mbox{ closed and convex}, \Leb(C) \le\delta
\}
< \frac{c}{c+1} ,
\]
whence
\[
L(\phi,Q) \to- \infty\qquad\mbox{as } \max_{x \in\R^d}
\phi
(x) \to\infty.
\]
Note also that $L(\phi,Q) \le\max_{x \in\R^d} \phi(x)$ for any
$\phi
\in\Phi(Q)$. These considerations show that $L(Q)$ is finite and, for
suitable constants $M_o < M_*$, equals the supremum of $L(\phi,Q)$ over
all $\phi\in\Phi(Q)$ such that $M_o \le\max_x \phi(x) \le M_*$.

Next we show the existence of a maximizer $\phi\in\Phi(Q)$ of
$L(\cdot
,Q)$. Let $(\phi_n)_n$ be a sequence of functions in $\Phi(Q)$ such
that $- \infty< L(\phi_n,Q) \uparrow L(Q)$ as $n \to\infty$, where
$M_n := \max_{x \in\R^d} \phi_n(x) \in[M_o,M_*]$ for all $n \ge1$.
Now we show that
%
%
\begin{equation}
\label{ineq:uniform lower bound}
\inf_{n \ge1} \phi_n(x_o) > - \infty\qquad
\mbox{for any } x_o \in\interior(\csupp(Q)) .
\end{equation}
If $\phi_n(x_o) < M_n$, then $x_o$ is not an interior point of the
closed, convex set $\{\phi_n \ge\phi_n(x_o)\}$. Hence
\begin{eqnarray*}
\int\phi_n \,dQ
& \le& \phi_n(x_o) + \bigl(M_n - \phi_n(x_o)\bigr) Q \{\phi_n \ge\phi
_n(x_o)\}
\\
& \le& \phi_n(x_o) + \bigl(M_n - \phi_n(x_o)\bigr) h(Q,x_o) \\
& \le& \phi_n(x_o)\bigl(1 - h(Q,x_o)\bigr) + \max(M_n,0)
\end{eqnarray*}
with $h(Q,x_o) < 1$ defined in Lemma \ref{lem:weak semicontinuity of
csupp}. In the case of $\phi_n(x_o) = M_n$ these inequalities are true
as well. Thus
\[
\phi_n(x_o) \ge- \frac{\max(M_n,0) - L(\phi_n,Q)}{1 - h(Q,x_o)}
\ge- \frac{\max(M_*,0) - L(\phi_1,Q)}{1 - h(Q,x_o)} ,
\]
which establishes (\ref{ineq:uniform lower bound}). Combining (\ref
{ineq:uniform lower bound}) with $\phi_n \le M_*$, we may deduce from
Lem\-ma~3.3 of
\citet{Schuhmacheretal2009}
that there exist constants $a$ and $b > 0$ such that
%
%
\begin{equation}
\label{ineq:uniform subexponentiality}
\phi_n(x) \le a - b \|x\|\qquad
\mbox{for all } n \in\mathbb{N}, x \in\R^d .
\end{equation}

The inequalities (\ref{ineq:uniform lower bound}) and (\ref{ineq:uniform
subexponentiality}) and Lemma \ref{lem:limits of sequences in GG} with
$C \supset\interior(\csupp(Q))$ and $\bar{\phi}(x) := a - b\|x\|$
imply existence of a function $\psi\in\Phi$ and a subsequence $(\phi
_{n(k)})_k$ of $(\phi_n)_n$ such that $\psi= - \infty$ on $\R^d
\setminus\csupp(Q)$ and
\begin{eqnarray*}
\limsup_{k \to\infty} \phi_{n(k)}(x)
& \le& \psi(x) \le a - b\|x\|\qquad
\mbox{for all } x \in\R^d , \\
\lim_{k \to\infty} \phi_{n(k)}(x)
& = & \psi(x) > - \infty
\qquad\mbox{for all } x \in\interior(\csupp(Q)) .
\end{eqnarray*}
Since the boundary of $\csupp(Q)$ has Lebesgue measure zero, it follows
from dominated convergence that $\int e^{\psi(x)} \,dx = 1$. Moreover,
applying Fatou's lem\-ma to the nonnegative functions $x \mapsto a - b\|
x\| - \phi_{n(k)}(x)$ yields
\[
\limsup_{k \to\infty} \int\phi_{n(k)} \,dQ
\le\int\psi \,dQ .
\]
Hence
\[
L(Q) \ge L(\psi,Q) \ge\limsup_{k \to\infty} L\bigl(\phi
_{n(k)},Q\bigr) = L(Q)
\]
and thus $L(\psi,Q) = L(Q)$.

Uniqueness of the maximizer $\psi$ follows essentially from strict
convexity of the exponential function: if $\tilde{\psi} \in\Phi(Q)$
with $L(\tilde{\psi},Q) > - \infty$, then $L( (1 - t)\psi+
t\tilde
{\psi},Q )$ is strictly concave in $t \in[0,1]$, unless $\Leb\{
\psi\ne\tilde{\psi}\} = 0$. But for $\psi, \tilde{\psi} \in\Phi(Q)$,
the latter requirement is equivalent to $\psi= \tilde{\psi}$ everywhere.
\end{pf*}

In our proofs of Theorems \ref{thm:characterization for d=1} and \ref
{thm:mallows-continuity} we utilize a special approximation scheme for
functions in $\Phi$:
\begin{Lemma}[(DSS 2010)]
\label{lem:approximation of phi from above}
For any function $\phi\in\Phi$ with nonempty domain and any
parameter $\eps> 0$ set
\[
\phi^{(\eps)}(x)
:= \inf_{(v,c)} (v^\top x + c)
\]
with the infimum taken over all $(v,c) \in\R^d \times\R$ such that
$\|
v\| \le\eps^{-1}$ and $\phi(y) \le v^\top y + c$ for all $y \in\R^d$.
This defines a function $\phi^{(\eps)} \in\Phi$ which is real-valued
and Lipschitz-continuous with constant $\eps^{-1}$. Moreover, it
satisfies $\phi^{(\eps)} \ge\phi$ with equality if and only if
$\phi
$ is real-valued and Lipschitz-continuous with constant $\varepsilon
^{-1}$. In general, $\phi^{(\eps)} \downarrow\phi$ pointwise as
$\eps
\downarrow0$.
\end{Lemma}
\begin{pf*}{Proof of Theorem \ref{thm:characterization for d=1}}
Let $P$ be the distribution corresponding to~$F$. Suppose first that
$\phi= \psi(\cdot| Q)$. Then it follows from (\ref
{ineq:characterization2}) and Fubini's theorem that
\begin{eqnarray*}
0
& = & \int_{\R} x (Q - P)(dx) \\
& = & \int_{\R} \int_{\R} ( 1 \{0 < t < x\} - 1\{x \le t \le
0\}
) \,dt
(Q - P)(dx) \\
& = & \int_{\R} \bigl( 1 \{0 < t\} (F - G)(t)
- 1\{t \le0\} (G - F)(t) \bigr) \,dt \\
& = & \int_{\R} (F - G)(t) \,dt .
\end{eqnarray*}
Moreover, for any $x \in\R$, the function $s \mapsto(s - x)^+$ is
convex so that (\ref{ineq:characterization1}) and Fubini's theorem yield
\[
0
\le\int_{\R} (s - x)^+ (Q - P)(ds)
= - \int_{-\infty}^x (F - G)(t) \,dt .
\]

It remains to be shown that $\int_{-\infty}^x (F - G)(t) \,dt \ge0$
for $x \in\SSS(\phi)$. Suppose first that $x \in\interior(\dom
(\phi
))$. Note that $\phi' := \phi'(\cdot +)$ is nonincreasing on the
interior of $\dom(\phi)$ with
\[
\phi(x_2) - \phi(x_1) = \int_{x_1}^{x_2} \phi'(u) \,du
\qquad\mbox{for } x_1, x_2 \in\interior(\dom(\phi)) \mbox
{ with } x_1 < x_2 .
\]
Moreover, $x \in\SSS(\phi)$ implies that $\phi'(x - \delta) > \phi
'(x +
\delta)$ for all $\delta> 0$ satisfying $x \pm\delta\in\interior
(\dom(\phi))$. For such $\delta> 0$ we define
\[
H_\delta(s) := \int_{-\infty}^s H_\delta'(u) \,du
\]
with
\[
H_\delta'(u)
:= \cases{
0,
&\quad \mbox{for} $u \le x-\delta$, \vspace*{2pt}\cr
\dfrac{\phi'(x - \delta) - \phi'(u)}{\phi'(x-\delta) - \phi
'(x+\delta)},
&\quad for $x-\delta< u \le x+\delta$, \vspace*{2pt}\cr
1, &\quad for $u \ge x+\delta$.}
\]
One can easily verify that $\phi+ t H_\delta$ is upper semicontinuous
and concave whenever $0 < t \le\phi'(x-\delta) - \phi'(x+\delta)$. In
case of $t < - \inf_{u \in\R} \phi'(u)$ it is also coercive. Thus it
follows from (\ref{ineq:characterization}) that
\begin{eqnarray*}
0 &\le&\int_{\R} H_\delta(s) (P - Q)(ds)
\to \int_{\R} (s - x)^+ (P - Q)(ds) \qquad(\delta\downarrow
0) \\
& = & \int_{-\infty}^x (F - G)(t) \,dt .
\end{eqnarray*}
When $x \in\SSS(\phi)$ is the left or right endpoint of $\dom(\phi)$,
we define $\Delta(s) := (s - x)^+$ and conclude analogously that $\int
_{-\infty}^x (F - G)(t) \,dt \ge0$.

Now suppose that the distribution function $F$ with log-density $\phi
\in\Phi$ satisfies the integral (in)equalities stated in Theorem \ref
{thm:characterization for d=1}. Let $\Delta\dvtx\R\to\R$ be
Lipschitz-continuous with constant $L$, so for arbitrary $x,y \in\R$
with $x < y$,
\[
\Delta(y) - \Delta(x) = \int_x^y \Delta'(t) \,dt
\]
with $\Delta' \dvtx\R\to[-L,L]$ measurable. Then
\[
\int\Delta d(Q - P)
= \int_\R\Delta'(t) (F - G)(t) \,dt .
\]
Since $\int(F - G)(t) \,dt = 0$, we may continue with
\begin{eqnarray*}
\int\Delta d(Q - P)
& = & \int_\R\bigl(\Delta'(t) + L\bigr) (F - G)(t) \,dt \\
& = & \int_\R\int_{-L}^L 1\{s < \Delta'(t)\} \,ds (F - G)(t)
\,dt \\
& = & \int_{-L}^L \int_{A(\Delta',s)} (F - G)(t) \,dt \,ds
\end{eqnarray*}
with $A(\Delta',s) := \{ t \in\R\dvtx\Delta'(t) > s\}$. Now we apply
this representation\break to the function $\Delta:= \phi^{(\eps)}$ for some
$\eps> 0$, that is, $L = \eps^{-1}$. Here one can show that $A(\Delta
',s)$ equals either $\varnothing$ or $\R$ or a half-line with right
endpoint $a(\phi,s) = \min\{t \in\R\dvtx\phi'(t +) \le s\}$. But this
entails that $a(\phi,s) \in\SSS(\phi)$, whence
$\int_{A(\Delta',s)} (F - G)(t) \,dt = 0$ for all $s \in(-L,L)$. Consequently,
\[
\int\phi^{(\eps)} \,d(Q - P) = 0 .
\]
If we consider $\Delta:= \psi^{(\eps)}$ with $\psi:= \psi(\cdot
|
Q)$, the sets $A(\Delta',s)$ are still half-lines with right endpoint
or empty or equal to $\R$. Thus
$\int_{A(\Delta',s)} (F - G)(t) \,dt \le0$ for all $s \in(-L,L)$, whence
\[
\int\psi^{(\eps)} \,d(Q - P) \le0 .
\]
Since $\phi^{(\eps)} \downarrow\phi$ and $\psi^{(\eps)}
\downarrow\psi
$ as $\eps\downarrow0$, and since $\int\phi \,dP$ and $\int\psi
\,dQ$ exist in $\R$, we can deduce from monotone convergence that $\int
\phi \,d(Q - P) = 0 \ge\int\psi \,d(Q - P)$. Since $\int e^{\phi
(x)} \,dx = \int e^{\psi(x)} \,dx = 1$, this entails that
\[
L(\phi,Q) = L(\phi,P)
\ge L(\psi,P) \ge L(\psi,Q) ,
\]
where the first displayed inequality follows from log-concavity of $P$
with log-density $\phi$. Thus $\phi= \psi$.
\end{pf*}

Theorem \ref{thm:weak upper semicontinuity} and the second part of
Theorem \ref{thm:mallows-continuity} are a consequence of the
following result:
\begin{Theorem}
\label{thm:semi- and full continuity}
Let $(Q_n)_n$ be a sequence of distributions in $\QQ_o$ such that $Q_n
\to_w Q \in\QQ_o$, $L(Q_n) \to\lambda\in[-\infty,\infty]$ and
$\int\|x\| Q_n(dx) \to\gamma\in[0,\infty]$ as $n \to\infty$.
Then $\gamma\ge\int\|x\| Q(dx)$, and $\lambda> -\infty$ if and
only if $\gamma< \infty$. Moreover,
\[
\lambda\cases{
< L(Q), &\quad if $\displaystyle \gamma> \int\|x\| Q(dx) $, \cr
= L(Q) \in\R, &\quad if $\displaystyle \gamma= \int\|x\| Q(dx) < \infty
$.}
\]
In the latter case, the densities $f := \exp\mbox{${}\circ{}$}\psi(\cdot| Q)$
and $f_n := \exp\mbox{${}\circ{}$}\psi(\cdot| Q_n)$ are well defined for
sufficiently large $n$ and satisfy
\begin{eqnarray*}
\lim_{n \to\infty, x \to y} f_n(x) & = & f(y)
\qquad\mbox{for all } y \in\R^d \setminus\partial\{f > 0\} , \\
\limsup_{n \to\infty, x \to y} f_n(x) & \le& f(y)
\qquad\mbox{for } y \in\partial\{f > 0\} , \\
\lim_{n \to\infty} \int| f_n(x) - f(x) | \,dx
& = & 0 .
\end{eqnarray*}
\end{Theorem}

Before presenting the proof of this result, let us recall two
elementary facts about weak convergence and unbounded functions:
\begin{Lemma}
\label{lem:weak convergence and unbounded f}
Suppose that $(Q_n)_n$ is a sequence in $\QQ$ converging weakly to some
distribution $Q$. If $h$ is a nonnegative and continuous function on~$\R^d$, then
\[
\liminf_{n \to\infty} \int h \,dQ_n \ge\int h \,dQ .
\]
If the stronger statement $\lim_{n \to\infty} \int h \,dQ_n = \int h
\,dQ < \infty$ holds, then
\[
\lim_{n \to\infty} \int f \,dQ_n = \int f \,dQ
\]
for any continuous function $f$ on $\R^d$ such that $|f|/(1 + h)$ is bounded.
\end{Lemma}
\begin{pf*}{Proof of Theorem \ref{thm:semi- and full continuity}}
The asserted inequality $\gamma\ge\int\|x\| Q(dx)$ follows from
the first part of Lemma \ref{lem:weak convergence and unbounded f} with
$h(x) := \|x\|$.

Suppose that $\gamma< \infty$. Then with $\phi(x) := - \|x\|$,
\[
\lambda\ge\lim_{n\to\infty} L(\phi,Q_n)
= - \gamma- \int e^{-\|x\|} \,dx + 1 > - \infty.
\]
In other words, $\lambda= - \infty$ entails that $\gamma= \infty$.

From now on suppose that $\lambda> - \infty$, and without loss of
generality let $L(Q_n) > - \infty$ for all $n \in\mathbb{N}$. We have
to show that $\gamma< \infty$ and that $\lambda\le L(Q)$ with
equality if and only if $\gamma= \int\|x\| Q(dx)$. To this end
we analyze the functions $\psi_n := \psi(\cdot| Q_n)$ and their
maxima $M_n := \max_{x \in\R^d} \psi_n(x)$. First of~all,
%
%
\begin{equation}
\label{eq:bounded maxima}
(M_n)_n \mbox{ is bounded} .
\end{equation}
This can be verified as follows: since $L(Q_n) = \int\psi_n \,dQ_n
\le M_n$, the sequence $(M_n)_n$ satisfies $\liminf_{n \to\infty} M_n
\ge\lambda$. With similar arguments as in the proof of Theorem \ref
{thm:existence and uniqueness} one can deduce that $(M_n)_n$ is bounded
from above, provided that
\[
\limsup_{n \to\infty} Q_n(C_n) < 1
\]
for any sequence of closed and convex sets $C_n \subset\R^d$ with
$\lim
_n \Leb(C_n) = 0$. To this end we refer to the proof of Lemma \ref
{lem:csupp} in [DSS 2010]: there exist a simplex $\tilde{\Delta} =
\conv
(\tilde{x}_0,\ldots,\tilde{x}_d)$ with positive Lebesgue measure and
open sets $U_0$, $U_1, \ldots, U_d$ with $Q(U_j) \ge\eta> 0$ for $0
\le j \le d$, such that $\tilde{\Delta} \subset C$ for any convex set
$C$ with $C \cap U_j \ne\varnothing$ for $0 \le j \le d$. But $\liminf_n
Q_n(U_j) \ge Q(U_j) \ge\eta$ for all $j$. Hence $\Leb(C_n) < \Leb
(\tilde{\Delta})$ entails that $Q_n(C_n) \le1 - \min_{0 \le j \le d}
Q_n(U_j) \le1 - \eta+ o(1)$ as $n \to\infty$.

Another key property of the functions $\psi_n$ is that
%
%
\begin{equation}
\label{eq:lower bounds}
\liminf_{n \to\infty} \psi_n(x_o) > - \infty
\qquad\mbox{for any } x_o \in\interior(\csupp(Q)) .
\end{equation}
For
\[
L(Q_n) = \int\psi_n \,dQ_n
\le\psi_n(x_o) + \bigl(M_n - \psi_n(x_o)\bigr) h(Q_n,x_o) ,
\]
whence as $n \to\infty$,
\[
\psi_n(x_o)
\ge- \frac{\max(M_n,0) - L(Q_n)}{1 - h(Q_n,x_o)}
\ge- \frac{\limsup_{\ell\to\infty} \max(M_\ell,0) -
\lambda
}{1 - h(Q,x_o)} + o(1)
\]
by virtue of Lemma \ref{lem:weak semicontinuity of csupp}.
Combining (\ref{ineq:uniform lower bound}) with (\ref{eq:bounded
maxima}) we may again deduce that there exist constants $a$ and $b > 0$
such that
%
%
\begin{equation}
\label{ineq:uniform subexponentiality 2}
\psi_n(x) \le a - b \|x\|\qquad
\mbox{for all } n \in\mathbb{N}, x \in\R^d .
\end{equation}

As in the\vspace*{1pt} proof of Theorem \ref{thm:existence and uniqueness} we can
replace $(Q_n)_n$ with a subsequence such that for suitable constants
$a$, $b > 0$ and a function $\tilde{\psi} \in\Phi$ the following
conditions are met: $\interior(\csupp(Q)) \subseteq\dom(\tilde
{\psi})$ and
\begin{eqnarray*}
\psi_n(y), \tilde{\psi}(y)
& \le& a - b \|y\| \qquad\mbox{for all } y \in\R^d, n \in\mathbb{N}
, \\
\lim_{n \to\infty, x \to y} \psi_n(x)
& = & \tilde{\psi}(y) \qquad\mbox{for all } y \in\interior(\dom
(\tilde
{\psi})), \\
\limsup_{n \to\infty, x \to y} \psi_n(x)
& \le& \tilde{\psi}(y) \qquad\mbox{for all } y \in\R^d .
\end{eqnarray*}
In particular,
\[
\lambda= \lim_{n\to\infty} \int\psi_n \,dQ_n
\le\lim_{n \to\infty} \int(a - b\|x\|) Q_n(dx)
= a - b \gamma,
\]
whence
\[
\gamma< \infty.
\]
Moreover, $\int\exp(\tilde{\psi}(x)) \,dx = \lim_{n \to\infty}
\int
\exp(\psi_n(x)) \,dx = 1$, by dominated convergence.

By Skorohod's theorem, there exists a probability space $(\Omega
,\mathcal{A},\mathbb{P})$ with random variables $X_n \sim Q_n$ and $X
\sim Q$ such that $\lim_{n \to\infty} X_n = X$ almost surely. Hence
Fatou's lemma, applied to the random variables $H_n := a - b \|X_n\| -
\psi_n(X_n)$, yields
\begin{eqnarray*}
\lambda &=& \lim_{n \to\infty} \int\psi_n \,dQ_n
= \lim_{n \to\infty}
\biggl( \int(a - b\|x\|) \,dQ_n - \Ex(H_n) \biggr) \\
& \le& a - b \gamma- \Ex\Bigl( \liminf_{n \to\infty} H_n
\Bigr)
\\ & \le& a - b\gamma- \Ex\bigl( a - b\|X\| - \tilde{\psi}(X)
\bigr) \\
& = & b \biggl( \int\|x\| Q(dx) - \gamma\biggr) + \int\tilde
{\psi
}(x) Q(dx) \\
& \le& b \biggl( \int\|x\| Q(dx) - \gamma\biggr) + L(Q) .
\end{eqnarray*}
Thus $\lambda< L(Q)$ if $\gamma> \int\|x\| Q(dx)$.

It remains\vspace*{1pt} to analyze the case $\gamma= \int\|x\| Q(dx) < \infty$.
Here $\lambda\le L(\tilde{\psi},Q) \le L(Q)$, and it remains to show
that $\lambda\ge L(Q)$ which would entail that $\tilde{\psi}$ equals
the unique maximizer $\psi:= \psi(\cdot| Q)$. With the
approximations $\psi^{(1)} \ge\psi^{(\eps)} \ge\psi$, $0 < \eps
\le
1$, introduced in Lemma \ref{lem:approximation of phi from above}, it
follows from their Lipschitz-continuity and Lemma \ref{lem:weak
convergence and unbounded f} that $\lambda= \lim_{n \to\infty}
L(\psi
_n,Q_n)$ is not smaller than
\[
\lim_{n \to\infty} L\bigl(\psi^{(\eps)},Q_n\bigr)
= L\bigl(\psi^{(\eps)},Q\bigr) = \int\psi^{(\eps)} \,dQ - \int\exp
\bigl(\psi
^{(\eps)}(x)\bigr) \,dx + 1 .
\]
By monotone convergence, applied to the functions $\psi^{(1)} - \psi
^{(\eps)}$, and dominated convergence, applied to $\exp\mbox{${}\circ{}$} \psi
^{(\eps)}$,
\[
\lambda\ge\lim_{\eps\downarrow0} L\bigl(\psi^{(\eps)},Q\bigr) = L(\psi,Q) =
L(Q) .
\]

Note that the probability densities $f = \exp\mbox{${}\circ{}$}\psi$ and $f_n =
\exp\mbox{${}\circ{}$}\psi_n$ obviously satisfy
\begin{eqnarray*}
\lim_{n \to\infty, x \to y} f_n(x)
& = & f(y) \qquad\mbox{for all } y \in\R^d \setminus\partial\{f >
0\}, \\
\limsup_{n \to\infty, x \to y} f_n(x)
& \le& f(y) \qquad\mbox{for all } y \in\partial\{f > 0\} .
\end{eqnarray*}
In particular, $(f_n)_n$ converges to $f$ almost everywhere w.r.t.
Lebesgue measure, whence $\int| f_n(x) - f(x) | \,dx \to0$.

The only problem is that we established these properties only for a
\textit{subsequence} of the original sequence $(Q_n)_n$. But elementary
considerations outlined in [DSS 2010] show that this is sufficient.
\end{pf*}
\begin{pf*}{Proof of Theorem \ref{thm:mallows-continuity}}
The assertions of this theorem are essentially covered by Theorem \ref
{thm:semi- and full continuity} as long as $Q \in\QQ_o\cap\QQ^1$. It
only remains to show that $L(Q_n) \to\infty$ if $D_1(Q_n,Q) \to0$ for
some $Q \in\QQ^1 \setminus\QQ_o$. Thus $\int\|x\| Q(dx) <
\infty$
and $Q(H) = 1$ for a hyperplane $H = \{x \in\R^d \dvtx u^\top x = r\}$
with a unit vector $u \in\R^d$ and some $r \in\R$. For $k \ge1$ we
define $\phi_k \in\Phi$ via
\[
\phi_k(x) := - \|a_k + B_k x\| + \log(k) ,
\]
where $B_k := I - uu^\top+ k uu^\top$ is a real, $d \times d$ matrix
and $a_k := - kr u$. Note that $\det(B_k) = k$ and $\phi_k(x) = \log(k)
- \|x\|$ for $x \in H$. Thus
\[
L(\phi_k,Q_n)
\to L(\phi_k,Q)
= \log(k) - \int\|x\| Q(dx) + \int e_{}^{-\|x\|} \,dx .
\]
Since the right-hand side may be arbitrarily large, $\lim_{n\to\infty}
L(Q_n) = \infty$.
\end{pf*}
\begin{pf*}{Proof of Theorem \ref{thm:existence regression}}
Note that $\bs{v} \mapsto\hat{Q}_{\bs{v}}$ defines a continuous
mapping from $\R^n$ into the space of probability distributions on $\R$
with finite first moment, equipped with Mallows distance $D_1$.
Moreover, by our assumption that $\bs{Y} \notin\MM(\bs{x})$, none of
the distributions $\hat{Q}_{m(\bs{x})}$, $m \in\MM$, degenerates to a~Dirac
measure. According to Theorem \ref{thm:mallows-continuity}, the
mapping $\bs{v} \mapsto L(\hat{Q}_{\bs{v}})$ is thus continuous from
$\MM(\bs{x})$ into~$\R$.

When proving existence of a maximizer, as explained in Section \ref
{ssec:maximum likelihood estimation}, we may restrict our attention to
the closed subset $\MM(\bs{x},\bar{Y}) := \{ \bs{v} \in\MM
(\bs
{x}) \dvtx\bar{v} = \bar{Y} \}$ of $\MM(\bs{x})$, where generally
$\bar{w}$ denotes the arithmetic mean $n^{-1} \sum_{i=1}^n w_i$ for a
vector $\bs{w} \in\R^n$. But for $\bs{v} \in\MM(\bs{x},\bar{Y})$,
\[
\int| x - \mu(\hat{Q}_{\bs{v}}) | \hat{Q}_{\bs{v}}(dx)
= \frac{1}{n} \sum_{i=1}^n |Y_i - v_i|
\ge\frac{1}{n} \sum_{i=1}^n |v_i| - \frac{1}{n} \sum_{i=1}^n
|Y_i| ,
\]
and the right-hand side tends to infinity as $\|\bs{v}\| \to\infty$.
Thus it follows from Lem\-ma~\ref{lem:life is simple} that
\[
L(\hat{Q}_{\bs{v}}) \to- \infty\qquad
\mbox{as } \|\bs{v}\| \to\infty, \bs{v} \in\MM(\bs{x},\bar
{Y}) ,
\]
and this coercivity, combined with continuity of $\bs{v} \mapsto
L(\hat
{Q}_{\bs{v}})$ and $\MM(\bs{x},\bar{Y})$ being closed, yields the
existence of a maximizer.
\end{pf*}
\begin{pf*}{Proof of Theorem \ref{thm:convolution}}
The proof that $Q \star R \in\QQ_o \cap\QQ^1$ is elementary and
omitted here. By affine equivariance (Remark \ref{rem:equivariance}),
we may and do assume that $\int y R(dy) = 0$. Now let $\psi:= \psi
(\cdot| Q)$ and $\tilde{\psi} := \psi(\cdot| Q\star R)$. Then
\[
L(Q\star R)
= \int\!\!\int\tilde{\psi}(x+y) Q(dx) R(dy) \\
= \int\tilde{\psi}_R \,dQ ,
\]
where
\[
\tilde{\psi}_R(x) := \int\tilde{\psi}(x + y) R(dy)
\le\tilde{\psi}(x)
\]
by Jensen's inequality. Hence
\[
L(Q\star R) \le\int\tilde{\psi} \,dQ = L(\tilde{\psi},Q)
\le L(Q) .
\]

Now suppose that $L(Q\star R) = L(Q)$, so in particular, $\tilde{\psi}
= \psi$. It follows from $\tilde{\psi}_R \le\tilde{\psi} \in\Phi
$ and
Fatou's lemma that $\tilde{\psi}_R \in\Phi$ with $\int\exp(\tilde
{\psi
}_R(x)) \,dx \le\int\exp(\tilde{\psi}(x)) \,dx = 1$. Thus
\[
L(Q)
= L(Q\star R)
\le L(\tilde{\psi}_R,Q)
\le L(Q) ,
\]
that is, $\tilde{\psi}_R = \psi= \tilde{\psi}$ and
%
%
\begin{equation}
\label{eq:convolution psi}
\psi(x) = \int\psi(x + y) R(dy)\qquad
\mbox{for all } x \in\R^d .
\end{equation}
It remains to be shown that (\ref{eq:convolution psi}) entails $R =
\delta_0$.
Note that $K := \{ x \in\R^d \dvtx\psi(x) = M_o \}$ with $M_o
:= \max_{y \in\R^d} \psi(y)$ defines a compact set. Hence for any unit
vector $u \in\R^d$ there exists a vector $x(u) \in K$ such that
$u^\top x(u) \ge u^\top x$ for all $x \in K$. But then $\psi(x(u) + y)
< M_o$ for all $y \in\R^d$ with $u^\top y > 0$. Hence
\[
M_o = \psi(x(u)) = \int\psi\bigl(x(u) + y\bigr) R(dy)
\]
implies that $R\{y \dvtx u^\top y > 0\} = 0$. Since $u$ is an arbitrary
unit vector, this entails that $\csupp(R) = \{0\}$, that is, $R =
\delta_0$.
\end{pf*}
\begin{pf*}{Proof of Theorem \ref{thm:consistency regression}}
Assumptions (A.2) and (A.3) imply that the empirical distribution $\hat{Q}_n :=
\hat{Q}_{n,\mu_n}$ of the true errors $\eps_{ni}$ satisfies both
$D_{\mathrm{BL}}(\hat{Q}_n,Q) \to_p 0$ and $\int|t| \hat{Q}_n(dt) \to_p
\int
|t| Q(dt)$. Thus $D_1(\hat{Q}_n,Q) \to_p 0$.

To verify the assertions of the theorem it suffices to consider a
sequence of \textit{fixed} vectors $\bolds{\eps}_n = (\eps_{ni})_{i=1}^n
\in
\R^n$ such that for a constant $c > 0$ to be specified later,
%
%
\begin{equation}
\label{eq:nice eps_n}
D_1(\hat{Q}_n,Q) + \sup_{m \in\MM\dvtx \|(m - \mu_n)(\bs{x}_n)\|_n
\le c}
D_{\mathrm{BL}} ( \hat{Q}_{n,m}, Q_{n,m} )
\to0 .
\end{equation}
Our goal is to show that $(\hat{f}_n,\hat{\mu}_n)$, viewed as a
function of $\bolds{\eps}_n$ and thus fixed, too, is well defined for
sufficiently large $n$ with
%
%
\begin{equation}
\label{eq:consistency eps_n}
\int| \hat{f}_n(x) - f(x) | \,dx \to0
\quad\mbox{and}\quad
\| (\hat{\mu}_n - \mu_n)(\bs{x}_n) \|_n \to0 .
\end{equation}
Note that we replaced $f_n$ with $f = \exp\mbox{${}\circ{}$} \psi(\cdot|
Q)$ because $\int|f_n(x) - f(x)| \,dx$ tends to $0$.

We know already that we have to restrict our attention to the set $\hat
{\MM}_n$ of all $m \in\MM_n$ such that $\int t \hat{Q}_{n,m}(dt) =
0$, that is, $\int t R_{(\mu_n - m)(\bs{x}_n)}(dt) = - \int t
\hat
{Q}_n(dt)$ converges to $0$. Since $\{m(\bs{x}_n) \dvtx m \in\hat{\MM
}_n\}
$ is a closed subset of $\R^n$ by (A.1), we may argue as in the proof
of Theorem \ref{thm:existence regression} that a maximizer $\hat{\mu
}_n$ of $L(\hat{Q}_{n,m})$ over all $m \in\hat{\MM}_n$ does exist. It
is possible that $L(\hat{Q}_{n,\hat{\mu}_n}) = \infty$, but if we can
show that $D_1(\hat{Q}_{n,\hat{\mu}_n},Q) \to0$, then $\hat{f}_n$
exists for sufficiently large $n$, too. Thus we may rephrase (\ref
{eq:consistency eps_n}) as
%
%
\begin{equation}
\label{eq:consistency eps_n 2}
D_1(\hat{M}_n,Q) \to0
\quad\mbox{and}\quad
\int|t| \hat{R}_n(dt) \to0 ,
\end{equation}
where $\hat{M}_n := \hat{Q}_{n,\hat{\mu}_n}$ and $\hat{R}_n :=
R_{(\mu
_n - \hat{\mu}_n)(\bs{x}_n)}$.

Note first that $\check{\mu}_n := \mu_n + \int t \hat{Q}_n(dt)$
belongs to $\hat{\MM}_n$, whence
%
%
\begin{equation}
\label{eq:LhM_n lower}
L(\hat{M}_n) \ge L(\hat{Q}_{n,\check{\mu}_n}) = L(\hat{Q}_n) \to
L(Q)
\end{equation}
by Theorem \ref{thm:mallows-continuity}. On the other hand
\[
\int|t| \hat{M}_n (dt)
= \frac{1}{n} \sum_{i=1}^n | \eps_{ni} + (\mu_n - \hat
{\mu
}_n)(x_{ni}) |
\ge\int|t| \hat{R}_n(dt) - \int|t| \hat{Q}_n(dt) .
\]
Thus, by Lemma \ref{lem:life is simple}, $\hat{\mu}_n$ satisfies
$\int|t| \hat{R}_n(dt) \le c$ for sufficiently large $n \in
\mathbb
{N}$, provided that
$c$ is larger than $\int|t| Q(dt) + \exp(- L(Q))$. In particular,
\[
D_{\mathrm{BL}}(\hat{M}_n, Q_n \star\hat{R}_n)
= D_{\mathrm{BL}}(\hat{M}_n, Q_{n,\hat{\mu}_n})
\to0 .
\]
Since $D_{\mathrm{BL}}(Q_n \star\hat{R}_n, Q\star\hat{R}_n) \le D_{\mathrm{BL}}(Q_n,Q)
\to0$, we know that even
\[
D_{\mathrm{BL}}(\hat{M}_n, Q \star\hat{R}_n)
\to0 .
\]

Since $(\hat{R}_n)_n$ is tight, to verify (\ref{eq:consistency eps_n
2}) we may consider a subsequence $(\hat{R}_{n(k)})_k$ that converges
weakly to some distribution $R$ as $k \to\infty$. Then
$\hat{M}_{n(k)} \to_w Q\star R$, so
\[
\limsup_{k\to\infty} L\bigl(\hat{M}_{n(k)}\bigr)
\le L(Q\star R) \le L(Q)
\]
by Theorems \ref{thm:weak upper semicontinuity} and \ref
{thm:convolution}. Because of (\ref{eq:LhM_n lower}) we even know that
$L(\hat{M}_{n(k)}) \to L(Q\star R) = L(Q)$ as $k \to\infty$.
Consequently, we may deduce from Theorems \ref{thm:weak upper
semicontinuity} and~\ref{thm:convolution} that
\[
\lim_{k\to\infty} D_1\bigl(\hat{M}_{n(k)}, Q\star R\bigr) = 0
\quad\mbox{and}\quad
R = \delta_a \qquad\mbox{for some } a \in\R.
\]

It remains to be shown that $a = 0$ and $\lim_{k\to\infty} \int|t|
\hat{R}_{n(k)}(dt) = 0$. Elementary arguments reveal that for arbitrary
$r > 0$ and $n \in\mathbb{N}$,
\begin{eqnarray*}
\int|t| \hat{M}_{n}(dt)
& \ge& \int\min(|t|,r) \hat{M}_{n}(dt)
+ \int|t| \hat{R}_{n}(dt) \\
&&{} - \int\min(|t|,2r) \hat{R}_{n}(dt)
- \int(|t| - r)^+ \hat{Q}_{n}(dt) .
\end{eqnarray*}
Hence $\int|t| \hat{R}_{n(k)}(dt)$ is not greater than
\begin{eqnarray*}
&&\int\min(|t|,2r) \hat{R}_{n(k)}(dt)
+ \int(|t| - r)^+ \hat{M}_{n(k)}(dt) + \int(|t| - r)^+ \hat
{Q}_{n(k)}(dt) \\
&&\qquad \to \int\min(|t|,2r) R(dt)
+ \int(|t| - r)^+ Q\star R(dt) + \int(|t| - r)^+ Q(dt)
\end{eqnarray*}
as $k \to\infty$. As $r \uparrow\infty$, the limit on the right-hand
side converges to $\int|t| R(dt) = |a|$. Consequently, $\lim_{k\to
\infty} D_1(\hat{R}_{n(k)},R) = 0$. But then $0 = \lim_{k\to\infty}
\int t \hat{R}_{n(k)}(dt)$ coincides with $\int t R(dt) = a$.
\end{pf*}

In our proofs of Theorems \ref{thm:consistency linear regression} and
\ref{thm:consistency isotonic regression} we utilize a simple
inequality for the bounded Lipschitz distance in terms of the
Kolmogorov--Smirnov distance,
\[
D_{\mathrm{KS}}(Q,Q') := \sup_{t \in\R} | (Q' - Q)((-\infty,t])
|,
\]
of two distributions $Q, Q' \in\QQ(1)$:
\begin{Lemma}[(DSS 2010)]
\label{lem:DBL and DKS}
Let $Q$ and $Q'$ be distributions on the real line. Then for arbitrary
$r > 0$,
\[
D_{\mathrm{BL}}(Q,Q') \le4 Q\bigl(\R\setminus(-r,r]\bigr) + 4(r+1) D_{\mathrm{KS}}(Q,Q') .
\]
\end{Lemma}
\begin{pf*}{Proof of Theorem \ref{thm:consistency linear regression}}
A key insight is that the empirical distributions $\hat{Q}_{n,m}$ are
close to their expectations $Q_{n,m}$ with respect to
Kolmogorov--Smirnov distance, uniformly over all $m \in\MM_n$. Namely,
\begin{eqnarray*}
&&{\sup_{m \in\MM_n, r \in\R}}
| (\hat{Q}_{n,m} - Q_{n,m})((-\infty,r]) | \\
&&\qquad = \sup_{b \in\R^{q(n)}, s \in\R}
\Biggl| \frac{1}{n} \sum_{i=1}^n \bigl( 1 \{Y_{ni} - b^\top x_{ni}
\le
s \}
- \Pr(Y_{ni} - b^\top x_{ni} \le s) \bigr) \Biggr| \\
&&\qquad \le {\sup_{H \in\HH_n} }| (\hat{M}_n - M_n)(H) | ,
\end{eqnarray*}
where $\HH_n$ denotes the family of all closed half-spaces in $\R
^{q(n)+1}$ while $\hat{M}_n$ is the empirical distribution of the
random vectors $(Y_{ni},x_{ni}^\top)_{}^\top\in\R^{q(n)+1}$, $1 \le i
\le n$, and $M_n := \Ex\hat{M}_n$. Now we utilize well-known results
from empirical process theory: $\HH_n$ is a Vapnik--\v{C}ervonenkis
class with VC-dimension $q(n)+3$, and $\hat{M}_n$ is the arithmetic
mean of $n$ independent random probability measures. Thus
\[
\Ex\sup_{m \in\MM_n} D_{\mathrm{KS}}(\hat{Q}_{n,m},Q_{n,m})
\le C \sqrt{ \frac{q(n)+3}{n} }
\]
for some universal constant $C$ [see
\citet{Pollard1990}, Theorems 2.2 and 3.5,
and
\citet{vdVaartWellner1996}, Theorem 2.6.4 and Lemma 2.6.16].

Since for fixed $c > 0$ the family $\{ Q_{n,m} \dvtx n \in\mathbb{N},
m \in\MM_n $ with $\|(m - \mu_n)\times(\bs{x}_n)\|_n \le c
\}$
is tight, the previous finding, combined with Lemma \ref{lem:DBL and
DKS}, implies that
\[
\lim_{n\to\infty} \Ex\sup_{m \in\MM_n \dvtx \|(m - \mu
_n)(\bs
{x}_n)\|_n \le c}
D_{\mathrm{BL}}(\hat{Q}_{n,m},Q_{n,m})
= 0 .
\]
\upqed
\end{pf*}

\section*{Acknowledgments}
Constructive comments by an Associate Editor and two referees are
gratefully acknowledged.


%
\printaddresses

\end{document}